\documentclass[12pt,twoside,reqno]{amsart}
\usepackage{amsxtra,amscd}
\usepackage{graphicx}
\usepackage{amsmath}
\usepackage{amsfonts}
\usepackage{amssymb}

\newtheorem{theorem}{Theorem}[section]

\theoremstyle{definition}
\newtheorem{definition}[theorem]{Definition}

\newtheorem{remark}[theorem]{Remark}

\numberwithin{equation}{section}

\begin{document}
\title{on the algebraic fundamental groups}
\author{Feng-Wen An}
\address{School of Mathematics and Statistics, Wuhan University, Wuhan,
Hubei 430072, People's Republic of China}
\email{fwan@amss.ac.cn}

\begin{abstract}
Passing from arithmetic schemes to algebraic schemes, in a similar manner we will
have the computation of the \'{e}tale fundamental group of an algebraic
scheme and then will define and discuss the qc fundamental group of an
algebraic scheme in this paper. The qc fundamental group will  also give a prior estimate of the \'{e}tale
fundamental group.
\end{abstract}

\maketitle

\begin{center}
{\tiny {Contents} }
\end{center}

{\tiny \qquad {Introduction} }

{\tiny \qquad {1. Preliminaries}}

{\tiny \qquad {2. The \'{E}tale Fundamental Group}}

{\tiny \qquad {3. The qc Fundamental Group} }

{\tiny \qquad {References}}

\section*{Introduction}

As an arithmetic scheme behaves frequently like the ring of algebraic
integers for a number field, it has been seen that there exist many tricks
arising from algebraic number theory which still work for the case of an
arithmetic scheme; hence, we have obtained some properties of the profinite
fundamental groups of arithmetic schemes (for example, see \cite{An2,An3,An4,An5,An6}).

Passing from arithmetic schemes to algebraic schemes, in this paper we will
have the computation of the \'{e}tale fundamental group of an algebraic
scheme; then we will define and discuss the qc fundamental group of an
algebraic scheme which will also give a prior estimate of the \'{e}tale
fundamental group. These results are related to the Section Conjecture of
Grothendieck (see \cite{faltings}).

\subsection*{Convention}

By an \textbf{algebraic variety} we will understand an
integral scheme $X$ over a field $k$ of finite type  in the paper.
In such a case, $X$ is also said to be an \textbf{algebraic $k-$variety}.
Here, the field $k$ can have an arbitrary characteristic.

\subsection*{Acknowledgment}

The author would like to express his sincere gratitude to Professor Li
Banghe for his advice and instructions on algebraic geometry and topology.

\bigskip

\section{Preliminaries}

For convenience, let us fix notation and definitions in this subsection.

\subsection{Notation}

Fixed an integral domain $D$. In the paper, we let $Fr(D)$ denote the field
of fractions on $D$.

If $D$ be a subring of a field $\Omega $, the field $Fr(D)$ will always
assumed to be contained in $\Omega $.

Let $E$ be an extension of a field $F$ (not necessarily algebraic). $E$ is
said to be \textbf{Galois} over $F$ if $F$ is the fixed subfield of the
Galois group $Gal(E/F)$.

\subsection{Affine Covering with Values}

Fixed a scheme $X$. As usual, an affine covering of the scheme $X$ is a
family
\begin{equation*}
\mathcal{C}_{X}=\{(U_{\alpha },\phi _{\alpha };A_{\alpha })\}_{\alpha \in
\Delta }
\end{equation*}
such that for each $\alpha \in \Delta $, $\phi _{\alpha }$ is an isomorphism
from an open set $U_{\alpha }$ of $X$ onto the spectrum $Spec{A_{\alpha }}$
of a commutative ring $A_{\alpha }$. Each $(U_{\alpha },\phi _{\alpha
};A_{\alpha })\in \mathcal{C}_{X}$ is called a \textbf{local chart}.

An affine covering $\mathcal{C}_{X}$ of $X$ is said to be \textbf{reduced}
if $U_{\alpha }\neq U_{\beta }$ holds for any $\alpha \neq \beta $ in $%
\Delta $.

Let $\mathfrak{Comm}$ be the category of commutative rings with identity.
Fixed a subcategory $\mathfrak{Comm}_{0}$ of $\mathfrak{Comm}$. An affine
covering $\{(U_{\alpha },\phi _{\alpha };A_{\alpha })\}_{\alpha \in \Delta }$
of $X$ is said to be \textbf{with values} in $\mathfrak{Comm}_{0}$ if 
 for each $\alpha \in \Delta $ there are $\mathcal{O}_{X}(U_{\alpha})=A_{\alpha}$ and $U_{\alpha}=Spec(A_{\alpha})$, where
 $A_{\alpha }$ is a ring contained in $\mathfrak{Comm}_{0}$.

In particular, let $\Omega $ be a field and let $\mathfrak{Comm}(\Omega )$
be the category consisting of the subrings of $\Omega $ and their
isomorphisms. An affine covering $\mathcal{C}_{X}$ of $X$ with values in $%
\mathfrak{Comm}(\Omega )$ is said to be \textbf{with values in the field $%
\Omega $}.

Assume that $\mathcal{O}_{X}$ and $\mathcal{O}^{\prime}_{X}$ are two structure sheaves on the underlying space of an integral scheme $X$. The two integral schemes $(X,\mathcal{O}_{X})$ and $(X, \mathcal{O}^{\prime}_{X})$ are said to be \textbf{essentially equal} provided that for any open set $U$ in $X$, we have
 $$U \text{ is affine open in }(X,\mathcal{O}_{X}) \Longleftrightarrow \text{ so is }U \text{ in }(X,\mathcal{O}^{\prime}_{X})$$ and in such a case,  $D_{1}=D_{2}$ holds or  there is $Fr(D_{1})=Fr(D_{2})$ such that for any nonzero $x\in Fr(D_{1})$, either $$x\in D_{1}\bigcap D_{2}$$ or $$x\in D_{1}\setminus D_{2} \Longleftrightarrow x^{-1}\in D_{2}\setminus D_{1}$$ holds, where $D_{1}=\mathcal{O}_{X} (U)$ and $D_{2}=\mathcal{O}^{\prime}_{X} (U)$.

 Two schemes $(X,\mathcal{O}_{X})$ and $(Z,\mathcal{O}_{Z})$ are said to be \textbf{essentially equal} if the underlying spaces of $X$ and $Z$ are equal and the schemes $(X,\mathcal{O}_{X})$ and $(X,\mathcal{O}_{Z})$ are essentially equal.

\subsection{Quasi-Galois Closed}

Fixed a field $k$. Let $X$ and $Y$ be algebraic $k-$varieties and let $f:X\rightarrow Y$ be a
surjective morphism of finite type. Denote by $Aut\left( X/Y\right) $ the group of
automorphisms of $X$ over $Z$.

By a \textbf{conjugate} $Z$ of $X$ over $Y$, we understand an algebraic $k-$variety $Z$ that is isomorphic to $X$ over $Y$.

\begin{definition}
$X$ is said to be \textbf{quasi-galois closed}  over $Y$
by $f$ if  there is an algebraically closed field $\Omega$
and a reduced affine covering $\mathcal{C}_{X}$ of $X$ with values in $
\Omega $ such that for any conjugate $Z$ of
$X$ over $Y$ the two conditions are satisfied:
\begin{itemize}
\item $(X,\mathcal{O}_{X})$ and $(Z,\mathcal{O}_{Z})$ are essentially equal if $Z$ has a reduced
affine covering with values in $\Omega$.

\item $\mathcal{C}_{Z}\subseteq \mathcal{C}_{X}$ holds if $\mathcal{C}_{Z}$
is a reduced affine covering of $Z$ with values in $\Omega $.
\end{itemize}
\end{definition}

\begin{remark}
We can prove the existence and the main property for an algebraic $k-$variety in an evident manner
(See \cite{An2,An3,An4,An5,An6}).
\end{remark}

\section{The \'{E}tale Fundamental Group}

\subsection{Definitions}

Fixed an algebraic variety $X$ over a field $k$. Let $L_{1}$ and $ L_{2}$ be two algebraic extensions over  $k(X)$, respectively.

\begin{definition}
$L_{2}$ is said to be a \textbf{finite $X$-formally
unramified Galois extension}  over $L_{1}$ if there are two algebraic $k$-varieties
$X_{1}$ and $X_{2}$ and a surjective morphism $f:X_{2}\rightarrow X_{1}$
such that
\begin{itemize}
\item $Sp[X]=Sp[X_{1}]=Sp[X_{2}]$, i.e., $X$,$X_{1}$, and $X_{2}$ have a same $sp$-completion.

\item $k\left( X_{1}\right) =L_{1}, \, k\left( X_{2}\right) =L_{2}$;

\item $X_{2}$ is a finite \'{e}tale Galois cover of $X_{1}$ by $f$.
\end{itemize}
In such a case, $X_{2}/X_{1}$ are said to be a \textbf{$X$-geometric model} of the field extension $L_{2}/L_{1}$.
\end{definition}

For $L=k(X)$, set

\begin{itemize}
\item $L^{al}\triangleq $ an algebraical closure of $L$;

\item $L^{sep}\triangleq $ the separable closure of $L$ contained in $L^{al}$%
;

\item $L^{au}\triangleq $ the union of all the finite $X$-formally unramified
subextensions over $L$ contained in $L^{al}$.
\end{itemize}

\begin{remark}
Let $L$ be a finitely generated extension over a number field $K$. Then $%
L^{au}$ is a subfield of $L^{al}$. In particular, it is seen that $L^{au}$
is a Galois extension over $L$.

Moreover, let $\omega \in L^{sep}$ be unramified over $L$. Then $f(\omega)\in
L^{sep}$ is also unramified over $L$ for any element $f$ of the absolute
Galois group $Gal(L^{al}/L)$.

Hence, by set inclusion, $L^{au}$ is (equal to and then defined to be) the
\textbf{maximal unramified subextensions} over $L$ (contained in $L^{al}$).
\end{remark}

\subsection{The Etale Fundamental Group}

By a trick similar to \cite{An4,An6}, we have the following result.

\begin{theorem}
Fixed any algebraic $k-$variety $X$. Then there exists an isomorphism
\begin{equation*}
\pi _{1}^{et}\left( X,s\right) \cong Gal\left( {k(X)}^{au}/k\left( X\right)
\right)
\end{equation*}%
between groups for any geometric point $s$ of $X$ over the separable closure
of the function field $k\left( X\right) .$
\end{theorem}

\section{The qc Fundamental Group}

\subsection{Definitions }

Let $X$ be an algebraic $k-$variety. Let $\Omega$ be a separably closed field  containing the function field $k\left( X\right) $. Here, $\Omega $ is not necessarily algebraic over $k\left(
X\right) .$

Define $X_{qc}\left[ \Omega \right] $ to be the set of algebraic $k-$varieties $%
Z$ satisfying the following conditions:
\begin{itemize}
\item $Z$ has a reduced
affine covering with values in $\Omega $;

\item there is a
surjective morphism $f:Z\rightarrow X$ of finite type such that $Z$ is
quasi-galois closed over $X.$
\end{itemize}

Set a partial order $\leq$ in the set $X_{qc}\left[ \Omega \right] $ in such
a manner:

Take any $Z_{1},Z_{2}\in X_{qc}\left[ \Omega \right] ,$ we say
\begin{equation*}
Z_{1}\leq Z_{2}
\end{equation*}
if there is a surjective morphism $\varphi :Z_{2}\rightarrow Z_{1}$ of
finite type such that $Z_{2}$ is quasi-galois closed over $Z_{1}.$

It is seen that $X_{qc}\left[ \Omega %
\right] $ is a directed set and
\begin{equation*}
\{Aut\left( Z/X\right) :Z\in X_{qc}\left[ \Omega \right] \}
\end{equation*}%
is an inverse system of groups. Hence, we have the following definition.

\begin{definition}
Let $X$ be an algebraic $k-$variety. Take any separably
closed field $\Omega $ containing $k\left( X\right) .$ The
inverse limit
\begin{equation*}
\pi _{1}^{qc}\left( X;\Omega \right) \triangleq {\lim_{\longleftarrow}}
_{Z\in X_{qc}\left[ \Omega \right] }{Aut\left( Z/X\right)}
\end{equation*}
of the inverse system $\{Aut\left( Z/X\right) :Z\in X_{qc}\left[ \Omega %
\right] \}$ of groups is said to be the \textbf{qc fundamental group} of the
scheme $X$ with coefficient in $\Omega .$
\end{definition}

\subsection{The qc Fundamental Group}

By a trick similar to \cite{An5}, we have the following results.

\begin{theorem}
Let $X$ be an algebraic $k-$variety. Take any separably
closed field $\Omega $ containing $k\left( X\right) .$ There are the
following statements.

$\left( i\right) $ There is a group isomorphism
\begin{equation*}
\pi _{1}^{qc}\left( X;\Omega \right) \cong Gal\left( {\Omega }/k\left(
X\right) \right) .
\end{equation*}

$\left( ii\right) $ Take any geometric point $s$ of $X$ over $\Omega $. Then
there is a group isomorphism
\begin{equation*}
\pi _{1}^{et}\left( X;s\right) \cong \pi _{1}^{qc}\left( X;\Omega \right)
_{et}
\end{equation*}%
where $\pi _{1}^{qc}\left( X;\Omega \right) _{et}$ is a subgroup of $\pi
_{1}^{qc}\left( X;\Omega \right) $. Moreover, $\pi _{1}^{qc}\left( X;\Omega
\right) _{et}$ is a normal subgroup of $\pi _{1}^{qc}\left( X;\Omega \right)
$.\bigskip
\end{theorem}

\begin{remark}
Let $X$ be an algebraic $k-$variety. Put
\begin{equation*}
\pi _{1}^{qc}\left( X \right)=\pi _{1}^{qc}\left( X;{k(X)}^{sep} \right).
\end{equation*}
Then there is a group isomorphism
\begin{equation*}
\pi _{1}^{qc}\left( X \right) \cong Gal\left( {k(X)}^{sep}/k\left( X\right)
\right).
\end{equation*}
\end{remark}

\begin{definition}
Let $X$ be an algebraic $k-$variety. The quotient group
\begin{equation*}
\pi _{1}^{br}\left( X\right) =\pi _{1}^{qc}\left( X;k\left( X\right)
^{sep}\right) /\pi _{1}^{qc}\left( X;k\left( X\right) ^{sep}\right) _{et}
\end{equation*}
is said to be the \textbf{branched group} of the algebraic variety $X$.
\end{definition}

The branched group $\pi _{1}^{br}\left( X\right)$ can reflect the
topological properties of the scheme $X,$ especially the properties of the
associated complex space $X^{an}$ of $X,$ for example, the branched covers
of $X^{an}$.

\begin{remark}
Let $X$ be an algebraic $k-$variety. Then we have
\begin{equation*}
\pi _{1}^{br}\left( X\right)=\{0\}
\end{equation*}
if and only if $X$ has no finite branched cover.
\end{remark}

\end{document}